\documentclass[12pt,leqno]{amsart}
\usepackage{graphicx}
\usepackage[utf8]{inputenc}
\usepackage[
backend=biber,
style=alphabetic,
sorting=nyt
]{biblatex}
\usepackage[T1]{fontenc}
\usepackage{lmodern}
\usepackage{geometry}
\usepackage{booktabs}
\usepackage{xcolor}
\usepackage{amsthm}
\usepackage{amsmath,amssymb}
\usepackage{listings}
\usepackage{emptypage}
\usepackage{newlfont}
\usepackage{enumerate}
\usepackage{mathrsfs}
\usepackage{hyperref}
\usepackage{tikz-cd}
\usepackage{comment}
\usepackage{cleveref}

\newcommand{\R}{\mathbb{R}}

\newcommand{\nexpL}{\exp^\bot_L}
\newcommand{\nexpS}{\exp^\bot_S}
\newcommand{\C}{\mathcal{C}}
\newcommand{\D}{\mathcal{D}}

\theoremstyle{plain} 
\newtheorem{thm}{Theorem}[section]
 
\newtheorem{lem}[thm]{Lemma} 
\newtheorem{prop}[thm]{Proposition}

\newtheorem*{thm*}{Theorem}
\newtheorem*{cor*}{Corollary} 
\newtheorem*{lem*}{Lemma} 
\newtheorem*{prop*}{Proposition}

\theoremstyle{definition}

\newtheorem*{defn*}{Definition}
\newtheorem*{ex*}{Example}

\theoremstyle{remark}

\allowdisplaybreaks[1]
\title[A generalization of the Gauss Lemma]
{A generalization of the Gauss Lemma}
\author{Lorenzo Scoffone}
\address{Karlsruher Institut f\"ur Technologie, Fakult\"at f\"ur Mathematik,
Englerstr.\ 2, 76131 Karlsruhe, Germany}
\email{lorenzo.scoffone@kit.edu}
\date{}
\hypersetup{
  pdftitle={A generalization of the Gauss Lemma},
  pdfauthor={Lorenzo Scoffone}
}

\begin{document}

\begin{abstract}
We give a generalization of the Gauss Lemma for normal exponential maps and a covering criterion as a consequence.
\end{abstract}

\maketitle

\section{Introduction}\label{sec:introduction}

The Gauss Lemma (cf \cite[Lemma 3.5]{doCarmo}) relates radial directions of the exponential map to lengths. This papers offers a generalization concerning rectifiable curves in the normal bundle of a $\C^{1,1}$ submanifold. Here and in the following, a $\C^1$ map is called $\C^{1,1}$ if its differentials depend locally Lipschitz continuously on the point. If $L$ is a submanifold of a Riemannian manifold $M$, let $\nexpL$ be its normal exponential map and $\pi:\nu L\to L$ its bundle map. Furthermore, consider the restriction of the norm of $TM$ to $\nu L$, $|\cdot|:\nu L\to[0,+\infty)$, $v\mapsto\|v\|_{T_{\pi(v)}M}$. Then, we prove the following:

\begin{lem}\label{lem:gauss1}
Consider a $\C^{1,1}$ submanifold $L$ of a Riemannian manifold $(M,g)$ and a rectifiable curve $\gamma:[0,a]\to\nu L$. Then $\ell(|\gamma|)\leq\ell(\nexpL\circ\gamma)$.
\end{lem}

This generalization has an application in extending a theorem by Hermann. In \cite{Hermann}, he proved that if the normal exponential map of a smooth submanifold of a smooth Riemannian manifold has no focal points, then it is a covering\footnote{In this paper \textit{covering} is intended in the classical sense, i.e. a continuous map $f:X\to Y$ such that any $y\in Y$ has an open neighbourhood $U$ such that $f^{-1}(U)$ is a disjoint union of some $V_i$, where each $f\vert_{V_i}:V_i\to U$ is a homeomorphism.} onto its image. This has been used in several contexts, especially in the theory of Riemannian foliations (for instance \cite{Hebda}). Motivated by the study of submetries (in which the best regularity of a general manifold fibre one can hope for is $\C^{1,1}$, see \cite{KapovitchLytchakSubmetries}), we obtain the following, replacing the differential requirement for the normal exponential map for a topological one:

\begin{thm}\label{thm:covering}
Assume there exists a closed $\C^1$ submanifold $L$ of a complete smooth Riemannian manifold $M$ such that $\nexpL:\nu L\to M$ is a local homeomorphism. Then $\nexpL$ is a covering.
\end{thm}

\subsection{Acknowledgments}
The author would like to thank Alexander Lytchak for suggesting the problem and for many helpful discussions.

\subsection{AI statement}
The author used OpenAI models for mathematical discussion and help with the exposition of the manuscript. The author takes full responsibility for the mathematical content of the manuscript.

\section{Proof of Lemma \ref{lem:gauss1}}\label{sec:gauss1}
Rectifiable is understood to be with respect to the metric that the smooth Riemannian manifold $TM$ induces on $\nu L$.
\begin{proof}[Proof of Lemma \ref{lem:gauss1}]
Parametrize $\gamma$ by arclength, so that it is Lipschitz. Then, by the Rademacher theorem $\gamma'(s_0)$ exists for almost all $s_0\in[0,a]$. It is enough to show that for any such $s_0$, we have
\begin{equation}\label{eq}
\left|\left.\frac{d}{ds}\right|_{s=s_0}|\gamma(s)|\right|\leq\left|\left.\frac{d}{ds}\right|_{s=s_0}(\nexpL\circ\gamma)(s)\right|.
\end{equation}
Note that \eqref{eq} depends only on $\gamma(s_0)$ and $\gamma'(s_0)$. Thus, fix $s_0$; we can ``smoothen'' $\gamma$ by considering a smooth $\widetilde\gamma:[0,a]\to TM$ such that $\gamma(s_0)=\widetilde\gamma(s_0)$ and $\gamma'(s_0)=\widetilde\gamma'(s_0)$.

\noindent
Choose a small $\varepsilon$ and for all $s\in(s_0-\varepsilon,s_0+\varepsilon)$ set $\alpha_s:[0,1]\to M,\ t\mapsto\exp(t\widetilde\gamma(s))$; define the Jacobi field $J_{s_0}(t)=\left.\frac{d}{ds}\right|_{s=s_0}\alpha_s(t)$. Note that the right hand side of \eqref{eq} equals $|J_{s_0}(1)|$. Furthermore,
\[
    J_{s_0}(0)=\left.\frac{d}{ds}\right|_{s=s_0}\pi(\widetilde\gamma(s))=d\pi_{\widetilde\gamma(s_0)}\widetilde\gamma'(s_0)=d\pi_{\gamma(s_0)}\gamma'(s_0)=(\pi\circ\gamma)'(s_0)\text{,}
\]
and is thus tangent to $L$. Therefore, it is orthogonal to $\gamma(s_0)=\alpha_{s_0}'(0)$. Note that this is the only time we are using that $\operatorname{Im}(\gamma)\subset\nu L$.

\noindent
Therefore we have the following:
\begin{align*}
\left|\left.\frac{d}{ds}\right|_{s=s_0}|\gamma(s)|\right|&=\left|\left.\frac{d}{ds}\right|_{s=s_0}\sqrt{g(\alpha_s'(0),\alpha_s'(0))}\right|=\left|\frac{\left.\frac{d}{ds}\right|_{s=s_0}g(\alpha_s'(0),\alpha_s'(0))}{2|\gamma(s_0)|}\right|\\
&=\left|\frac{g\left(\left.\frac{d}{ds}\right|_{s=s_0}\left.\frac{d}{dt}\right|_{t=0}\alpha_s(t),\alpha_{s_0}'(0)\right)}{|\gamma(s_0)|}\right|
=\left|\frac{g(J_{s_0}'(0),\gamma(s_0))}{|\gamma(s_0)|}\right|
\end{align*}
Now, decompose $J_{s_0}$ in its normal and parallel parts, namely $J_{s_0}(t)=J_{s_0}^n(t)+J_{s_0}^p(t)$. For the normal part we always have $J_{s_0}^n(0)\perp\alpha_{s_0}'(0)$ and $(J_{s_0}^n)'(0)\perp\alpha_{s_0}'(0)$. The first condition, together with $J_{s_0}(0)\perp\alpha_{s_0}'(0)$, implies then $J_{s_0}^p(0)=0$, while the second implies $g(J_{s_0}'(0),\alpha_{s_0}'(0))=g((J_{s_0}^p)'(0),\alpha_{s_0}'(0))$. Furthermore, for the parallel part $(J_{s_0}^p)'(0)$ is a multiple of $\alpha_{s_0}'(0)=\gamma(s_0)$, say $(J_{s_0}^p)'(0)=c\gamma(s_0)$ for some $c\in\R$, and $J_{s_0}^p(t)=J_{s_0}^p(0)+t(J_{s_0}^p)'(0)$, which in our case yields $J_{s_0}^p(1)=(J_{s_0}^p)'(0)$. Therefore,
\begin{align*}
\left|\left.\frac{d}{ds}\right|_{s=s_0}|\gamma(s)|\right|&=\frac{|g((J_{s_0}^p)'(0),\gamma(s_0))|}{|\gamma(s_0)|}=\frac{|c|\,|\gamma(s_0)|^2}{|\gamma(s_0)|}=|c|\,|\gamma(s_0)|\\
&=\left|(J_{s_0}^p)'(0)\right|=|J_{s_0}^p(1)|.
\end{align*}
But since the decomposition $J_{s_0}(1)=J_{s_0}^n(1)+J_{s_0}^p(1)$ is orthogonal, $\left|J_{s_0}^p(1)\right|\leq|J_{s_0}(1)|$, thus finishing the proof.
\end{proof}

\section{From $\C^1$ to $\C^{1,1}$}\label{sec:C1,1}

Approaching the proof of Theorem \ref{thm:covering}, let us see that our local homeomorphism assumption on $\nexpL$ alone guarantees more regularity for $L$ and $\nexpL$:
\begin{lem}\label{lem:C1,1}
Given a $\C^1$ submanifold $L$ of a Riemannian manifold $M$ such that $\nexpL:\D\to M$ is a local homeomorphism on its domain of definition $\D$, then $L$ is $\C^{1,1}$. In particular, $\nexpL$ is locally Lipschitz.
\end{lem}
\begin{proof}
Let us see that $L$ has positive reach: this implies that $L$ is
$\C^{1,1}$ cf. \cite[Proposition 6.1]{KapovitchLytchakSubmetries}.

\noindent
For any $x$ in $L$, since the topology of the bundle $\nu L$ is locally given by the product topology of $L$ and $\nu_xL=\R^k$ we can find a neighbourhood $U$ of $x$ in $L$ and an $r>0$ such that $\nexpL$ restricted to the tube $\nu^{\leq r}U=\{v\in\nu L\mid\pi(v)\in U,\ |v|\leq r\}$ is a homeomorphism. Call $V_r(x)=\nexpL(\nu^{\leq r}U)$.

\noindent
Assume by contradiction that there exist some $x_j\in M$ converging to $x$ with non-unique foot point projections to $L$, say $p_j^1$ and $p_j^2$. Then, clearly the $p_j^i$ converge to $x$, $i=1,2$. Call $v_j^i\in\nu_{p_j^i}L$ the normal vectors such that $x_j=\nexpL(v_j^i)$, $i=1,2$. Given $N>0$ such that $x_n\in V_r(x)$ for all $n\geq N$, we would have $v_N^i\in\nu^{\leq r}U$ and $x_N=\nexpL(v_N^1)\neq\nexpL(v_N^2)=x_N$. Thus, the foot point projection $\pi_L$ is well defined on an open  neighbourhood of $x$. Calling $O$ the union of all such neighbourhoods as $x$ varies in $L$, we therefore get that $\pi_L:O\to L$ is well defined: $L$ has positive reach.
\end{proof}

\section{A distance identity}\label{sec:dist}
The key point for the proof of Theorem \ref{thm:covering} is the following distance identity:
\begin{prop}\label{prop:dist}
Consider a $\C^1$ submanifold $S$ of a smooth Riemannian manifold $N$ such that $\nexpS:\D'\to N$ is a homeomorphism on its domain of definition $\D'\subset\nu S$. Then $d_S(\nexpS(v))=|v|$ for every $v\in\D'$.
\end{prop}

We first need to prove the following Lemma, which can again be seen as a generalization of the Gauss Lemma. In particular, it extends Lemma \ref{lem:gauss1} to the case of a priori non rectifiable curves when the normal exponential map is a \textit{global} homeomorphism.
\begin{lem}\label{lem:gauss2}
Consider a $\C^1$ submanifold $S$ of a Riemannian manifold $N$ such that $\nexpS:\D'\to N$ is a homeomorphism on its domain of definition $\D'\subset\nu S$. Then, for every curve $\eta:[0,a]\to\D'$, we have $\ell(|\eta\vert)\leq\ell(\nexpS\circ\eta)$.
\end{lem}

Note that $S$ and $N$ as in Lemma \ref{lem:gauss2} verify the hypotheses of both Lemma \ref{lem:gauss1} and Lemma \ref{lem:C1,1}. In particular, $S$ is $\C^{1,1}$, $\nu S$ is a locally Lipschitz submanifold of $TN$ and $\nexpS$ is locally Lipschitz. 

\begin{proof}[Proof of Lemma~\ref{lem:gauss2}]
Let $F$ be the local representation of $\nexpS$ in bi-Lipschitz coordinates for $\D'$ and smooth coordinates of $N$. Thus, $F$ is a homeomorphism between open sets of $\R^n$ that is Lipschitz, but possibly not bi-Lipschitz. Let us first prove the statement for $\eta$ when $\nexpS\circ\eta$ is a coordinate segment.

\noindent
The case $n=1$ is obvious, so assume without loss of generality that $n\geq2$. Consider a segment $[a,b]\subset\operatorname{Im}(F)$ and a small closed cylinder around it of radius $\delta$, denoted by $C$, so that also $C\subset\operatorname{Im}(F)$. $K:=F^{-1}(C)$ is compact since $F^{-1}$ is continuous.

\noindent
Let $W:=\mathrm{Span}\{b-a\}^\perp\cong\R^{n-1}$ be the orthogonal complement to the segment and $P:\R^n\to W$ be the projection onto it. The coarea formula (cf.\ \cite[Theorem 9.4]{Ambrosio}) gives
\[
\int_W\mathcal H^1\bigl(K\cap(P\circ F)^{-1}(y)\bigr)\,dy=\int_K J_{n-1}(P\circ F)(u)\,du<\infty,
\]
where $\mathcal H^1$ is the normalized one-dimensional Hausdorff measure and $J_{n-1}$ the $(n-1)$-Jacobian.

\noindent
Now, for $w\in W$ with $|w|\leq\delta$, consider the interval obtained by translating $[a,b]$ by $w$, namely $\{a+w+t(b-a)\mid t\in[0,1]\}$. Consider further its preimage through $F$, $\gamma_w(t)$ for $t\in[0,1]$. Note $\gamma_w([0,1])=F^{-1}(C\cap P^{-1}(P(a)+w))=K\cap(P\circ F)^{-1}(P(a)+w)$. From the coarea formula $\mathcal H^1(\gamma_w([0,1]))<\infty$ for almost every such $w$. Since $\gamma_w$ is continuous and injective, it is rectifiable for almost every $w$ by \cite[Theorem 2.6.2]{BBI}. Thus, its lift $\eta_w$ through the Lipschitz parametrization of $\D'$ is rectifiable for almost every $w$ as well. From Lemma~\ref{lem:gauss1} we then get $\ell(|\eta_w|)\leq\ell(\nexpS\circ\eta_w)$ for almost every $w$. Choosing such $w_j$ converging to $0$, smoothness of coordinates of $N$ implies that $\ell(\nexpS\circ\eta_{w_j})\to\ell(\nexpS\circ\eta_0)$. By lower semi-continuity of $\ell$ we get $\ell(|\eta_0|)\leq\liminf_j\ell(|\eta_{w_j}|)\leq\lim_j\ell(\nexpS\circ\eta_{w_j})=\ell(\nexpS\circ\eta_0)$. Thus we have proven the desired statement for coordinate segments.

\noindent
Finally, consider any curve $\eta:[0,a]\to\D'$. Assume $\ell(\nexpS\circ\eta)<\infty$ as otherwise there is nothing to prove. Consider a partition $0=t_0\leq t_1\leq\ldots\leq t_N=a$ such that each $c_i:=\nexpS\circ\eta\vert_{[t_i,t_{i+1}]}$ lies in a coordinate neighbourhood of $N$ like above. Each $c_i$ can be approximated by interpolation with coordinate segments $c_{ij}$ with the same endpoints. If the interpolation is fine enough, all segments lie in the same coordinate neighbourhood, and the sum of their lengths converges to the length of $c_i$. Their preimages through $\nexpS$ also converge uniformly to $\eta\vert_{[t_i,t_{i+1}]}$, and again by lower semi-continuity of $\ell$ we obtain $\ell(|\eta\vert_{[t_i,t_{i+1}]}|)\leq\ell(c_i)$ like in the above. Taking the sum over $i$ yields the desired statement.
\end{proof}

\begin{proof}[Proof of Proposition \ref{prop:dist}]
Consider any $v\in\D'$ and $x=\nexpS(v)$. Since $|\cdot|$ vanishes on the zero section of $S$, by Lemma \ref{lem:gauss2} we have $|v|\leq\ell(c)$ for every curve $c$ joining a point of $S$ to $x$. Thus $|v|\leq d_S(x)$. The converse is obvious considering the normal geodesic $t\mapsto\nexpS(tv)$ for $0\leq t\leq1$.
\end{proof}

\section{Proof of Theorem \ref{thm:covering}}\label{sec:covering}
We now consider on $\nu L$ the length structure induced by the local homeomorphism $\nexpL$: namely, given a curve $\widetilde\gamma$ on $\nu L$, we set $\ell(\widetilde\gamma)=\ell(\nexpL\circ\widetilde\gamma)$. Then, using the distance $\bar d$ associated to this length structure, $\nexpL$ is a local isometry; in particular, $\nu L$ with the pull-back metric through the local isometry is a smooth Riemannian manifold. Note however that as a submanifold of $TM$ it is only $\C^0$ (actually, by Lemma~\ref{lem:C1,1} it is $\C^{0,1}$). Note that the topology induced by $\bar d$ on $\nu L$ coincides with the original one on $\nu L$: since $\nexpL$ is a local homeomorphism, on any subset on which it is a homeomorphism a subset is open in both topologies if and only if its image is open. 

\noindent
Since we will be dealing with both the distance $d$ on $M$ and $\bar d$ on $N$, for this section we will adopt the notation $B_d(p,R)$ for metric balls.

\begin{proof}[Proof of Theorem \ref{thm:covering}]
Let us first see that $N:=\nu L$ and $S:=0_L$ verify the hypotheses of Proposition \ref{prop:dist}. For this section we will adopt for the sake of clarity the notation for normal bundles $\nu^M L$, $\nu^N S=\nu^{\nu^ML}0_L$.

\noindent
$0_L$ is obviously closed in $\nu^M L$. Furthermore, since $L$ is $\C^1$ and $\nexpL$ is a local isometry, $0_L$ is $\C^1$ as a submanifold of $\nu^M L$.

\noindent
For $v\in\nu_p^M L$, set $A(v)=(d(\nexpL)_{0_p})^{-1}v\in\nu_{0_p}^N S$. This defines a norm-preserving bundle homeomorphism $A:\nu^M L\to\nu^N S$ since $\nexpL$ is a local isometry and $d(\nexpL)_{0_p}(T_{0_p}S)=T_pL$. Consider the radial geodesic $t\mapsto tv$ in $N$, whose initial velocity is $A(v)$, so that $\nexpS(A(v))=v$. Since every normal vector to $S$ can be written as $A(v)$ and its corresponding geodesic is defined on $[0,1]$, the domain of definition of $\nexpS$ is all of $\nu^N S$ and $\nexpS=A^{-1}$ is a global homeomorphism onto $N$. 

\noindent
Thus we can apply Proposition \ref{prop:dist} to $0_L\subset\nu^ML$ to obtain $\bar d_{0_L}(v)=|A(v)|=|v|$ for all $v\in\nu^M L$.

\medskip
We can now prove that $(\nu^M L,\bar d)$ is complete. Fix $0_p\in0_L$; if $v\in\bar B_{\bar d}(0_p,R)$, then by the above $|v|=\bar d_S(v)\leq\bar d(v,0_p)\leq R$. Furthermore, $d(\nexpL(v),p)\leq\bar d(v,0_p)\leq R$. Now, $\nexpL(tv)$ for $0\leq t\leq1$ joins $\pi(v)$ to $\nexpL(v)$ and has length $|v|$. Therefore, 
\[
d(\pi(v),p)\leq d(\pi(v),\nexpL(v))+d(\nexpL(v),p)\leq |v|+R\leq2R.
\]
It follows that $\bar B_{\bar d}(0_p,R)\subset\{v\in\nu^M L: \pi(v)\in L\cap\bar B_d(p,2R),\ |v|\leq R\}$. Since the right hand side is compact, so is $\bar B_{\bar d}(0_p,R)$. 

\noindent
Thus given any $v\in\nu^M L$ and $R>0$, we get $\bar B_{\bar d}(v,R)\subset\bar B_{\bar d}(0_{\pi(v)},R+|v|)$, which is compact. Thus so is $\bar B_{\bar d}(v,R)$. Therefore, by the Hopf-Rinow-Cohn-Vossen Theorem for metric spaces \cite[Theorem 2.5.28]{BBI} $(\nu^M L,\bar d)$ is complete.

\medskip
Finally, note that $\nexpL$ is surjective: connect any $p\in M$ with a point $q\in L$ realizing $d(L,p)$, then the geodesic $\gamma$ from $q$ to $p$ is such that $\gamma'(0)\in\nu_q^M L$ and $\nexpL(\gamma'(0))=p$. By \cite[Theorem 3.4.18]{BBI} $\nexpL$ is then a covering.
\end{proof}

\printbibliography[heading=bibintoc, title={Bibliography}]

\end{document}